\documentclass[12pt,twoside]{article}
\usepackage[all]{xy}
\usepackage{a4,amsmath,amssymb,amsfonts,amscd,mathrsfs}
\reversemarginpar
\newcommand{\kkk}[1]{}

\newcounter{Abschnitt}[section]
\newcommand{\neu}[1]{\protect\refstepcounter{Abschnitt}\protect
   \label{#1}\vspace{1ex}
   {\bf \protect\arabic{section}.\protect\arabic{Abschnitt}}
                     \kkk{#1}}

\newcommand{\Alb}{{\rm Alb}}
\newcommand{\alb}{{\rm alb}}
\newcommand{\Aut}{{\rm Aut}}

\newcommand{\rk}{{\rm rk}}
\newcommand{\Hom}{{\rm Hom}}

\newcommand{\Pic}{{\rm Pic}}

\newcommand{\Ext}{{\rm Ext}}
\newcommand{\GL}{{\rm GL}}

\newcommand{\Spec}{{\rm Spec}}

\newcommand{\HNF}{{\rm HNF}}

\newcommand{\U}{{\rm U}}

\newcommand{\NE}{{\rm NE}}
\newcommand{\NS}{{\rm NS}}
\newcommand{\id}{{\rm id}}

\newcommand{\gr}{{\rm gr}}

\newcommand{\ch}{{\rm ch}}
\newcommand{\Td}{{\rm Td}}

\newcommand{\sh}{{\rm sh}}
\newcommand{\Sh}{{\rm Sh}}

\newcommand{\Ecal}{{\cal E}}

\newcommand{\Lcal}{{\cal L}}

\newcommand{\Ocal}{{\cal O}}

\newcommand{\Scal}{{\cal S}}

\newcommand{\cdop}{{\mathbb C}}

\newcommand{\qdop}{{\mathbb Q}}
\newcommand{\ndop}{{\mathbb N}}
\newcommand{\rdop}{{\mathbb R}}
\newcommand{\pdop}{{\mathbb P}}

\newcommand{\zdop}{{\mathbb Z}}

\newcommand{\shom}{{{\cal H} om}}

\newcommand{\dual}{^\lor}

\newcommand{\rarpa}[1]{\stackrel{#1}{\rightarrow}}

\newcommand{\Lrar}{\Leftrightarrow}

\newcommand{\proof}{{\bf Proof: }}

\newcommand{\qed}{{ \hfill $\Box$}}

\newcommand{\tip}{{\vspace{0.5em} }}
\author{Georg Hein}
\date{September 27, 2004}
\title{Generalized Albanese morphisms}

\begin{document}

\maketitle

\begin{abstract}
We define generalizations of the Albanese variety
for a projective variety $X$.
The generalized Albanese morphisms
$\xymatrix{\alb_r: X\ar@{-->}[r] & \Alb_r(X)}$
contract those curves 
$C$ in $X$ for which the induced morphism
$\Hom(\pi_1(X),\U(r)) \to \Hom(\pi_1(C),\U(r))$
has a finite image up to conjugation.
Thus, they may be interpreted as an $\U(r)$-version of the
Shafarevich morphism.
It is shown that for an algebraic surface $X$, we have a 
regular morphism $\alb_r:X \to \Alb_r(X)$ with the above property. 
\end{abstract}

\section{Introduction}
The Shafarevich conjecture states that the universal cover
$\tilde X$ of a projective variety $X$ is holomorphically
convex.
This means that there exists a proper surjection
$\tilde X \to Y$ to a Stein space. 
This implies (see \cite{Kol}) the existence of a morphism
$\sh_X:X \to \Sh(X)$ with connected fibers with the property
that a curve $C$ in $X$ is contracted by $\sh_X$, iff
the image of $\pi_1(C)$ in $\pi_1(X)$ is finite.
The Shafarevich variety $\Sh(X)$ is the quotient of the Stein space
$Y$ by the fundamental group $\pi_1(X)$ of $X$.

The Shafarevich conjecture is very hard to tackle in
terms of algebraic geometry because the holomorphic map
$\tilde X \to X$ is not a regular morphism of
algebraic varieties unless the fundamental group of $X$
is finite.
It is, however, possible to deal with the problem of finding the
morphism $\sh_X$.
A first version of such a morphism was given by Koll\'ar
in \cite{Kol}.
For K\"ahler manifolds, Campana constructs  \cite{Cam}
a generic version
of the Shafarevich morphism. Katzarkov
(see \cite{Kat1} and \cite{Kat2}) describes a
representation theoretical description of $\sh_X$,
for one representation $\pi_1(X) \to \GL(r)$.
We offer a different approach to such a description
for all representations
by considering generalized Albanese varieties.

Two constructions for the Albanese variety exist.
One construction (\ref{ALB}.\ref{alb1}) uses global holomorphic one forms,
while another (\ref{ALB}.\ref{alb2})
defines the Albanese
variety as the dual of the Picard torus.
The intermediate Jacobians are a generalization of the first construction.
Here we try to show how to generalize the second construction.
In the sense of Koll\'ar \cite{Kol},  where the Albanese variety can be see
as the abelian or $\U(1)$-version of the Shafarevich morphism,
we intend to give an $\U(r)$-version of the Shafarevich morphism
that is a map
$\alb_r(X):X \to \Alb_r(X)$ having connected fibers,
such that
an irreducible curve $C$ in $X$ is contracted under $\alb_r(X)$,
iff the image of the composition
$\Hom(\pi_1(X),\U(r)) \to \Hom(\pi_1(C),\U(r))$
has a finite image up to conjugation.
If the image of $\pi_1(C)$ in $\pi_1(X)$ is finite,
then $C$ is contracted by the generalized 
Albanese morphisms.
Therefore, the generalized Albanese morphisms
would factor through the Shafarevich morphism.

In section \ref{ALB}, we review the constructions
and some properties of the Albanese variety.

The main result of section \ref{CO}
is the construction of the nef line bundle $\Lcal_r$ on $X$
whose properties are investigated in section \ref{PROP}.
The main property of $\Lcal_r$ is that
it has degree zero exactly on those curves $C$ in $X$ for which 
the morphism $\Hom(\pi_1(X),\U(r)) \to \Hom(\pi_1(C),\U(r))$ has a finite image
up to conjugation (see theorem \ref{PROP}.\ref{Prop8}).
Thus, $\Lcal_r$ behaves like the pull back of some ample line bundle from
$\Alb_r(X)$ to $X$.

In the fifth section, we use the nef reduction to obtain a rational map
$\xymatrix{\alb_r: X \ar@{-->}[r] & \Alb_r(X)}$ with the above
property for the general fiber.  
In the sixth section, we go over the construction for
the case of a projective algebraic surface $X$.
Here we are able to deduce the existence of  regular
morphisms $\xymatrix{\alb_r: X \ar[r] & \Alb_r(X)}$
for all integers $r \geq 1$ (see \ref{SURFACE}.\ref{sfmain}).
In the last section, we give an example of a surface with a
classical Albanese variety of dimension one
and generalized Albanese varieties of dimension two. 
Thus, this example shows that the generalized Albanese
morphism may ``reveal'' more
than the classical Albanese morphism does.

{\bf Notations.}
We work with schemes over the complex numbers.
Since we need the restriction of (semi)stable
vector bundles to curves,
we are required to use the concept of Mumford--Takemoto or slope stability
for vector bundles.

\section{Two constructions for the Albanese variety}\label{ALB}
\neu{alb1}{\bf The classical construction for the Albanese variety.\/}
Here we assume that $X$ is a connected K\"ahler manifold.
We define the Albanese variety $\Alb(X)$ to be the quotient 
$$\Alb(X):= H^0(X,\Omega^1_X)\dual / H_1(X,\zdop) \,.$$
If we choose a point $x_0 \in X(\cdop)$,
then we can define the Albanese morphism $\alb_X:X \to \Alb(X)$ by $x
\mapsto \int_{\gamma_x}$ where $\gamma_x$ is a path connecting $x_0$ with
$x$.

\neu{alb2}{\bf The $\Pic^0(\Pic^0)$-description of $\Alb(X)$.\/}
Let $X$ be a smooth variety over an algebraically closed field $k$.
We consider the Picard torus $\Pic^0(X)$,
i.e.,
the component of $\Pic(X)$ containing $\Ocal_X$.
Furthermore,
we consider a Poincar\'e bundle $L$ on $X \times \Pic^0(X)$.
This bundle is not unique.
To normalize it we choose a point $x_0 \in X(k)$.
If we require that  $L|_{\{x_0\} \times \Pic^0(X)} \cong
\Ocal_{\Pic^0(X)}$,
then the Poincar\'e bundle $L$ is uniquely determined.
If we consider $L$ as a family of line bundles on $\Pic^0(X)$ parametrized
by
$X$,
then we obtain a morphism from $X$ to the Picard torus of $\Pic^0(X)$:
$$\alb_X: X \to \Pic^0(\Pic^0(X)) =: \Alb(X)\,.$$

\neu{alb3}{\bf Both constructions coincide for smooth varieties over
$\Spec(\cdop)$.\/}
The use of the same notations in the above constructions is justified,
because both coincide for a smooth projective variety over $\Spec(\cdop)$.
This follows from the universal property of the Albanese variety and the
duality between the Albanese variety and the Picard torus (cf.~\cite{GH}).
The following proposition describes the fibers of
the Albanese morphism.
Since it uses both descriptions
(\ref{ALB}.\ref{alb1} and \ref{ALB}.\ref{alb2}),
we have to assume that $X/\Spec(\cdop)$ is a smooth
projective variety.

\neu{alb4}{\bf Proposition.
(Description of the fibers of the Albanese morphisms
$\alb_X: X \to \Alb(X)$,
for a projective complex manifold $X$ (cf.~II.6 in \cite{GH}).)\/}\\
{\em 
If $\iota: Z \to X$ is a connected cycle in $X$,
then the following
conditions are equivalent:\\
(i) $Z$ is contained in a fiber of the morphism $\alb_X$;

\tip
(ii) The image of $\iota_*:H_1(Z,\zdop) \to H_1(X,\zdop)$ is finite;

\tip
(iii) The pull back morphism
$\iota^*:H^0(X,\Omega^1_X) \to H^0(Z,\Omega_Z^1)$ is trivial;

\tip
(iv) Let $\rho:\pi_1(X) \to \U(1)$
be a representation of the fundamental group.
Then the restriction $\rho|_{\pi_1(Z)}$ has a finite image;

\tip
(v) If  $L$ is a line bundle on $X \times S$,
then the pull back $\iota^*L$ on $Z \times S$ is of the form
$L_1 \boxtimes L_2$,
for any Noetherian scheme $S$.}

\section{The line bundle $\Lcal_r$}\label{CO}
\neu{Co1}{\bf The setup.}
We fix a smooth projective variety $X$
of dimension $n$ with a very ample line bundle $\Ocal_X(H)$
and a positive integer $r$.
Furthermore, we choose a geometric point $x_0 \in X$.
Let $M_r=M_X(r,0,0,\ldots,0)$ be the moduli space of
$S$--equivalence classes of slope semistable 
rank $r$ bundles $E$ with trivial Chern classes in $H^*(X,\zdop)$.
If $E$ is a vector bundle parametrized by $M_r$,
then we write $[E]$ for the corresponding point in $M_r(\cdop)$.
By the theorem of Uhlenbeck and Yau (see \cite{UY}), $M_r$ parameterizes
flat vector bundles on $X$ or representations of $\pi_1(X)$ in $\U(r)$
modulo conjugation. 
This implies that for $[E] \in M_r$, the restriction $E|_C$
of $E$ to any curve $C \subset X$
is semistable.
Moreover, $M_r$ is a projective scheme provided that we pass to
$S$--equivalence classes of semistable
bundles. This means we identify any vector bundle $E$ in a short exact sequence
$0 \to E' \to E \to E'' \to 0$ of slope zero bundles with the direct sum $E' \oplus E''$.
Thereafter, we will use use the symbol $M_r$ (or $M_r(X)$) for the projective
moduli space of $S$--equivalence classes of slope semistable bundles on $X$.

\neu{Co2}{\bf The line bundle $\Ocal_M(D_H)$.}
Using the polarization $H$ on $X$, we can define a polarization $D_H$ on $M$.
We choose a faithfully flat morphism $\psi:\tilde M \to M$, such that 
we have a universal sheaf $\tilde E$ 
on $\tilde M \times X$. That means, for any point $\tilde m \in \tilde M$,
the sheaf $\tilde E_{\tilde m}:=\tilde E|_{\{ \tilde m \} \times X} $ is a sheaf which
belongs to the $S$--equivalence class given by $\psi(\tilde m)$.
The theory of Quot schemes gives the existence of such morphisms.
Let $C = H_1 \cap H_2 \cap \ldots \cap H_{n-1}$ be a complete intersection
of $n-1$ divisors $H_i \in |H|$.
Furthermore,
we take a rank two vector bundle $F$ on $C$ with $\det(F) \cong \omega_C$.
We consider the following morphisms:
$$\xymatrix{ \tilde M \times X \ar[rr]^\Psi \ar[dd]_{\tilde q} \ar[rd]^{\tilde p}
&& M \times X \ar[ld]_p \ar[dd]_q\\
& X \\
\tilde M \ar[rr]^\psi &&M}$$

On $\tilde M$, we define the line bundle $\Ocal(\tilde D_H)$ to be the determinant of
cohomology
$$\Ocal_{\tilde M}(\tilde D_H):=
\det (\tilde q_! (\tilde E \otimes \tilde p^* F))^{-1}\,.$$
This line bundle descends to $M$, i.e. there exists a line bundle $\Ocal_M(D_H)$ and
an isomorphism $\Ocal_{\tilde M}(\tilde D_H) \cong \psi^* \Ocal_M(D_H)$.
Furthermore, the line bundle $\Ocal_M(D_H)$ does not depend on the choice of
the morphism $\psi:\tilde M \to M$ (see \cite{DN}). 

\neu{Co2a}{\bf Lemma.}
{\em 
The first Chern class $\tilde D_H:=c_1(\Ocal_{\tilde M}(\tilde D_H))$ is given by
$$\tilde D_H=\tilde q_*((2c_2(\tilde E)-c_1^2(\tilde E)).\tilde p^*H^{n-1})\,.$$
In particular, we have $D_{aH} = a^{n-1} D_H$.}

\tip
{\bf Proof:}
This is straightforward computation using the Grothendieck--Riemann--Roch
formula
$$\tilde D_H =-\left[ \tilde q_*( \ch(\tilde E) \cdot \Td(\tilde q) \cdot 
\tilde p^* \ch(F))\right]_1\,,$$
and the equalities
$$\begin{array}{rclcl}
\ch(\tilde E) &=& r + c_1(\tilde E)+\frac{c_1^2(\tilde E)-2c_2(\tilde E)}{2}+ \ldots \\
\Td(\tilde q) &=& \tilde p^* \Td(X) &=& 1-\frac{\tilde p^* K_X}{2} + \ldots \\
\ch(F) &=& \ch(\Ocal_X \oplus \Ocal_X(K_X+(n-1)H))\cdot \ch(\Ocal_H)^{n-1}
&=& 2H^{n-1}+ K_X.H^{n-1}\,.\\
\end{array}$$
Whereas the first two equalities are standard, the last equality follows from
the adjunction formula and the fact that a vector bundle $F$ on a curve $C$ is determined in
the Grothendieck group $K(C)$ by its rank and determinant.
\qed

\neu{Co2b}{\bf Remark.}
The line bundle $\Ocal_{M_r}(D_H)$ is an example of a generalized Theta
divisor (see \cite{DN} and \cite{Hei}).
This means that if for a given vector bundle $F$ on $C$
we have a dense
open subset $U \subset M_r$ such that
$[E] \in U$ implies $H^1(E|_C \otimes F)=0$,
then there exists a section $s_F \in \Ocal_{M_r}(D_H)$ which vanishes
exactly at the points $\{ [E] \in M_r \, | \, h^1(E|_C \otimes F) >0 \}$.
In this case, the divisor $\Ocal_{M_r}(D_H)$ is effective which justifies
this notation. It follows that $F$ itself is a semistable bundle.
However, the existence of such a bundle $F$ is not clear for $r >2$.

\neu{Co2c}{\bf Lemma.}
{\em
The line bundle $\Ocal_{M_r}(a \cdot D_H)$ is base point free for $a \geq r^2$.}

\tip
\proof
Let $[E] \in M_r$ be a geometric point.
Then $E|_C$ is a semistable vector bundle on the curve $C$.
Popa shows in \cite{Popa} that if $a \geq r^2$, then
there exists a rank $2a$ vector
bundle $G$ on $C$ with $\det(G) \cong \omega_C^{\otimes a}$ with
the property $H^*(C,G \otimes E|_C) = 0$. Thus, the generalized Theta
divisor $\theta_G \in H^0(M_r, \det(\tilde q_!(\tilde E \otimes \tilde p^*G))^{-1})$
associated to $G$ does not pass through the points in $\psi^{-1}[E]$. Analogously to
the computations in the proof of lemma \ref{CO}.\ref{Co2a},
we see that the line bundle 
$\det(\tilde q_!(\tilde E \otimes \tilde p^*G))^{-1}$ is isomorphic to
$\Ocal_{M_r}(a \cdot D_H)$.
\qed

\neu{Co2d}{\bf Lemma.}
{\em
The line bundle $\Ocal_{M_r}(D_H)$ is ample.}

\tip
\proof
The proof uses the fact that on a projective variety of dimension at least
two, the vector bundles $E$ and $E'$ with the same Hilbert polynomial
are isomorphic, if and only if
their restrictions to a sufficiently big ample divisor $H$ are isomorphic.
This follows from the long exact sequence
$$\xymatrix{\Hom(E,E'(-H)) \ar@{=}[d] \ar[r] & \Hom(E,E') \ar[r]
& \Hom(E,E'|_H) \ar[r] \ar@{=}[d] & \Ext^1(E,E'(-H)) \ar@{=}[d] \\
H^0(E\dual \otimes E'(-H)) && \Hom(E|_H,E'|_H) & H^1(E\dual \otimes E'(-H)).}$$
If we have a bounded family of vector bundles as in the case of those 
parametrized by $M_r$, then we can choose a divisor $H$ such that for any
two bundles in this family, the cohomology groups on the left and right hand vanish.

The restriction theorem of Mehta and Ramanthan (see \cite{MR})
tells us that, for a semistable
vector bundle $E$ and $H$ big enough, the formation of graded objects
commutes with restriction to $H$. Thus, we obtain an embedding
$\xymatrix{M_r= M_r(X) \ar@{^(->}[r] & M_r(H)}$.
Repeating the argument, we end up with an embedding
$\xymatrix{M_r \ar@{^(->}[r] & M_r(C)}$ for a complete intersection curve $C$.
By lemma \ref{CO}.\ref{Co2a}, we may assume that this curve $C$ is the
curve we considered in the construction of $\Ocal(D_H)$.

By construction, $\Ocal(D_H)$ is the pull back of the generalized
Theta line bundle on $M_r(C)$. This line bundle is known to be ample
by the work of Drezet and Narasimhan (see \cite{DN}).
Thus, the lemma holds.
\qed

\neu{Co3}{\bf The families $\Ecal_{r,i}$.}
Let $M_r^{\rm red} = \cup_{i=1}^l M_{r,i}$ be the decomposition of
the reduced scheme underlying $M_r$ into its
irreducible components, and let $\tilde M_{r,i}$ be the normalization
of the component $M_{r,i}$. We have a morphism
$\alpha_i:\tilde M_{r,i} \to M_r$ and consider the globally
generated line bundle (by lemma \ref{CO}.\ref{Co2c}) $N_{r,i}:=
\alpha_i^* \Ocal(r^2 \cdot D_H)$.
Let $C_{r,i}$ be the intersection of $\dim(\tilde M_{r,i})-1$
general global sections of $N_{r,i}$. By Bertini's theorem, $C_{r,i}$ can be assumed
to be a smooth irreducible curve.
Thus, by Langton's theorem (see \cite{Langton}), we have a universal vector
bundle $\Ecal_{r,i}$ on $C_{r,i} \times X$.
If the universal vector bundle $\Ecal_M$ on $M_r \times X$ existed,
then $\Ecal_{r,i}$ would be the pull back of this bundle to $C_{r,i} \times X$.

\neu{Co4}{\bf The line bundle $\Lcal_{r,i}$.}
We consider the vector bundle $\Ecal_{r,i}$ on $C_{r,i} \times X$ and the morphisms
$$\xymatrix{C_{r,i} & C_{r,i} \times X \ar[l]_-q \ar[r]^-p & X \,.}$$
We define the line bundle $N_{r,i}$ on the curve $C_{r,i}$ by
$N_{r,i} : = \det (\Ecal_{r,i}|_{C_{r,i} \times \{ x_0 \}} ) \,.$
Let $G_{r,i}$ be a vector bundle on $C_{r,i}$ with $\rk(G_{r,i})=2r$, and
$\det(G_{r,i}) \cong \omega_{C_{r,i}}^r \otimes N_{r,i}^{-2}$.
Similar to the definition of $\Ocal_{M_r}(D_H)$ in \ref{CO}.\ref{Co2},
we define the line bundle by
$$ \Lcal_{r,i} := \det(p_! ( \Ecal_{r,i} \otimes q^*G_{r,i} ))^{-1} \,.$$

\neu{Co4a}{\bf Remark.}
Unfortunately, in contrast to $\Ocal_{M_r}(D_H)$, the line bundle
$\Lcal_{r,i}$ is not independent of the choices. We next give an example
for this dependence on the choice of the family $\Ecal_{r,i}$.

Let
$X$ be a curve, and $\Ecal_{r,i}$ be a family of degree zero vector bundles
on $X$ parametrized by $C_{r,i}$. For a point $c \in C_{r,i}$, we consider
the vector bundle $\Ecal_c := \Ecal_{r,i}|_{\{ c \} \times X}$.
Furthermore, we assume that $\Ecal_c$ is not stable.
Thus, we have a
short exact sequence
$0 \to \Ecal'_c \to \Ecal_c \to \Ecal_c'' \to 0$
of degree zero vector bundles on $X$.
Denote by $\Ecal'_{r,i}$ the kernel of the natural surjection
$\Ecal_{r,i} \to \Ecal_c''$. The families $\Ecal'_{r,i}$ and $\Ecal_{r,i}$
parameterize the same $S$--equivalence classes of vector bundles on $X$.
A straightforward computation shows that the resulting line bundles $\Lcal_{r,i}$
and $\Lcal'_{r,i}$ fulfill
$$\Lcal'_{r,i} = \det(\Ecal_c')^{2 \rk(\Ecal_c'')} \otimes
\det(\Ecal_c'')^{-2\rk(\Ecal_c')} \otimes \Lcal_{r,i} \,.$$
Thus, we can only hope that the numerical type of $\Lcal_{r,i}$ is well
defined. This is the case as we will see in the next section
(see Corollary \ref{PROP}.\ref{Prop4}).

\neu{Co4b} We end this section
by defining the line bundle $\Lcal_r$ on $X$ by
$$\Lcal_r := \bigotimes_{i=1}^l \Lcal_{r,i} \,.$$

\section{Properties of the line bundles $\Lcal_r$}\label{PROP}\kkk{PROP}
\neu{Prop1}{\bf Relations defined by nef line bundles.}
Let $\Lcal$ be a nef line bundle on a proper variety $X$.
This line bundle defines an equivalence relation $\sim_\Lcal$
on the geometric points of $X$ as follows:
$$ x \sim_\Lcal x' \quad \Lrar \quad \left\{
\begin{array}{l}
\mbox{There exists a closed curve } C \subset X \\
\mbox{with }  x\in C,\, x' \in C, \mbox{ and } \Lcal.C=0.
\end{array} \right\}$$
We define the relation $\preceq$ on nef line bundles by:
The condition $L_1 \preceq L_2$ holds
if for any curve $C \subset X$ the inequality $L_1.C >0$ implies $L_2.C>0$.

We write $L_1 \prec L_2$ if $L_1 \preceq L_2$ holds, and there exists a curve
$C \subset X$ with $L_1.C = 0$ and $L_2.C>0$.
If $L_1 \preceq L_2$ holds, then the relation $\sim_{L_2}$ is contained in $\sim_{L_1}$,
i.e. if $x \sim_{L_2} x'$, then we have $x \sim_{L_1} x'$.
Whenever both relations $L_1 \preceq L_2$ and $L_2 \preceq L_1$ hold,
we write $L_1 \sim L_2$. This means that the relations $\sim_{L_1}$ and
$\sim_{L_2}$ coincide.

If we have a chain $L_0 \prec L_1 \prec \ldots \prec L_k$ of nef line bundles,
then the we have strict inclusions
$(\NE(X) \cap L_k^\bot) \subset (\NE(X) \cap L_{k-1}^\bot) \subset \ldots
\subset (\NE(X) \cap L_0^\bot)$ in the cone of curves. 
Since the orthogonal complements $ L_i^\bot$ are linear subspace of $H^2(X,\rdop)$,
we deduce that $k \leq \rho(X)= \dim(\NE(X))$.

\neu{Prop2}{\bf Theorem.} {\em The line bundle $\Lcal_{r,i}$ is nef, i.e.
for any morphism $\iota: Y \to X$ of a smooth curve $Y$ to $X$,
the degree of the line bundle $\iota^*\Lcal_{r,i}$ is 
non negative.
Furthermore, if the degree of $\iota^*\Lcal_{r,i}$ equals zero,
then for all geometric points $y_1$ and $y_2$ of $Y$ there
exists an isomorphism
$\Ecal_{C_{r,i} \times \{y_1\}} \cong \Ecal_{C_{r,i} \times \{y_2\}}$,
and all the vector bundles on $Y$ parametrized by $C_{r,i}$
are S--equivalent.}

\tip
\proof
We divide the proof into steps.

{\em Step 1: Reduction to the case where $X$ is a smooth projective curve.}\\
Since all the vector bundles parametrized by $M_r$ restrict to semistable bundles
on every closed subscheme of $X$, the pull back $(\id \times \iota)^* \Ecal_{r,i}$
is a family of semistable rank $r$ vector bundles on $Y$ parametrized by $C_{r,i}$.
Since the determinant of cohomology commutes with base change, we may assume $X=Y$.
Thus, we consider the vector bundle $\Ecal_{r,i}$ on the surface $C_{r,i} \times X$
and the following morphisms to smooth curves
$$\xymatrix{C_{r,i} & C_{r,i} \times X \ar[l]_-q \ar[r]^-p & X \,.}$$

{\em Step 2: The line bundles
$L_1$ and $L_2:=\Lcal_{r,i}$.}\\
Let $F$ be vector bundle on $X$ with $\rk(F)=2 r$ and $\det(F)=\omega_X^r$.
For a point $x_0 \in X$, we set
$N_{r,i}:=\det(\Ecal_{r,i}|_{C_{r,i} \times \{x_0 \}})$.
Let $G$ be a vector bundle on $C_{r,i}$ of rank $2r$ with
$\det(G)=\omega_{C_{r,i}}^r \otimes N_{r,i}^{-2}$.
We set $L_1:=\det(q_!(\Ecal_{r,i} \otimes p^*F))^{-1}$, and
$L_2:=\det(p_!(\Ecal_{r,i} \otimes q^*G))^{-1}$.
The nefness property of $\Lcal_{r,i}$ is equivalent to $\deg(L_2) \geq 0$.

\tip
{\em Step 3: $\deg(L_1)=\deg(L_2)$.}\\
We use the Grothendieck--Hirzebruch--Riemann--Roch theorem to compute the degrees
of the line bundles $L_1$ and $L_2$. Let us fix the notations before doing so:
By $c_0$ and $x_0$, we denote two geometric points of $C_{r,i}$ and $X$.
We use $F_p$ and $F_q$ to name the fibers $p^{-1}(x_0)$ and $q^{-1}(c_0)$.
The genera of $C_{r,i}$ and $X$ we denote by $g_C$ and $g_X$.
Since we are only interested in the degrees, we may assume
$\omega_X = \Ocal_X((2g_X-2)x_0)$ and $\omega_{C_{r,i}} = \Ocal_X((2g_C-2)c_0)$.
For the same reason, we have $\ch(F)=2r+2r(g_X-1)x_0$ and
$\ch(G)=2r+(2r(g_C-1)-2(\int_{C_{r,i} \times X}(F_p.c_1(\Ecal_{r,i}))))c_0$.
Furthermore, let $\Td(C_{r,i}) = 1-(g_C-1) c_0$ and $\Td(X)= 1-(g_X-1)x_0$
be the (numerical) Todd classes.
$$\begin{array}{rcl}
\deg(L_1) & = & -\int_{C_{r,i}} \ch( q_! (\Ecal_{r,i} \otimes p^* F)) \\
&=&-\int_{C_{r,i} \times X} \ch(\Ecal_{r,i} \otimes p^* F) p^* \Td(X) \\ 
&=& -\int_{C_{r,i} \times X} \ch(\Ecal_{r,i}) p^* \ch(F) p^*\Td(X) \\
&=& -\int_{C_{r,i} \times X} \ch(\Ecal_{r,i}) p^* (\ch(F) \Td(X)) \\
& = & -\int_{C_{r,i} \times X} 
\left( r+c_1(\Ecal_{r,i}) + \frac{c_1^2(\Ecal_{r,i}) - 2 c_2(\Ecal_{r,i})}{2} \right)
p^*(2r)\\
& = & r \cdot \int_{C_{r,i} \times X} (2 c_2(\Ecal_{r,i}) - c_1^2(\Ecal_{r,i}))
\end{array}$$
Analogously, we obtain for the degree of $L_2$ that:
$$\begin{array}{rcl}
\deg(L_2) & = & -\int_{C_{r,i} \times X}
\left( r+c_1(\Ecal_{r,i}) + \frac{c_1^2(\Ecal_{r,i}) - 2 c_2(\Ecal_{r,i})}{2} \right)
q^*(2r -( 2\int_{C_{r,i} \times X}(F_p.c_1(\Ecal_{r,i})c_0 )))\\
& = & r \cdot \int_{C_{r,i} \times X} (2 c_2(\Ecal_{r,i}) - c_1^2(\Ecal_{r,i}))
+2\cdot \left( \int_{C_{r,i} \times X} F_p.c_1(\Ecal_{r,i}) \right) \cdot
\left( \int_{C_{r,i} \times X} F_q.c_1(\Ecal_{r,i}) \right) 
\end{array}$$

Since $C_{r,i}$ parameterizes a family of $X$ vector bundles of degree zero,
the intersection number $\int_{C_{r,i} \times X} F_q.c_1(\Ecal_{r,i})$ equals zero.
Thus, we end up with the claimed equality.

\tip
{\em Step 4: $L_1^{\otimes r}$ is globally generated.
Thus, $\deg(L_2)=\deg(L_1) \geq 0$.}\\
This is a direct consequence of Popa's (see \cite{Popa})
result about the base point freeness
of certain powers of the generalized Theta line bundles. Here the semistability
of all vector bundles parametrized by $C_{r,i}$ is needed.
Popa's result uses Le Potier's method (cf. \cite{LPot}) of
constructing global sections of the generalized Theta bundle.

To prove $\deg(L_1) \geq 0$, we need fewer premises.
Indeed, if at least one point of $C_{r,i}$ para\-meterizes a semistable
vector bundle on $X$,
then a power of $L_1$ has a nontrivial section $s$ which is what we need.
This argument is the same as in the proof of Lemma \ref{CO}.\ref{Co2c}.

If the degree of $L_1$ is positive, then we consider two points $c_1$, $c_2$ of
$C_{r,i}$ with the property that the section $s$ vanishes at $s_1$ but not at $s_2$.
It follows, that the $X$ vector bundles parametrized by $c_1$ and $c_2$ are not
$S$--equivalent.

\tip
From now on, we assume that the degrees of $L_1$ and $L_2$ are zero.
For the following steps see also the proof of Theorem I.4 in
Faltings article \cite{Fal} or for the simpler rank two case
see Theorem 3.4 in \cite{Hei}. 

\tip
{\em Step 5: For any two geometric points $P$ and $Q$ of $X$, the vector bundles
$\Ecal_{r,i}|_{C_{r,i} \times \{ P \}}$ and
$\Ecal_{r,i}|_{C_{r,i} \times \{ Q \}}$ are isomorphic.}\\
We show that, given a point $P \in X$, for almost all points $Q \in X$, we have an
isomorphism between $\Ecal_{r,i}|_{C_{r,i} \times \{ P \}}$ and
$\Ecal_{r,i}|_{C_{r,i} \times \{ Q \}}$.
From this statement, the assertion of step 5 follows immediately.
Let $\Ecal_{c_0}:= \Ecal_{r,i}|_{\{ c_0 \} \times X}$ be the semistable
vector bundle on $X$ parametrized by $c_0$.
Take a vector bundle $F_P$ on $X$ such that $H^*(X, F_P \otimes \Ecal_{c_0}) = 0$.
Such a bundle exists because by Le Potier's result in \cite{LPot}.
Moreover, it defines a global section in a power of $L_1$ which does
not vanish at $c_0$. Since $\deg(L_1)=0$, this section has an empty vanishing divisor.
This implies that $H^*(X,\Ecal_{r,i}|_{\{ c \} \times X})=0$ for all points
$c \in C_{r,i}$. This implies $R^* q_*(\Ecal_{r,i} \otimes p^* F_P)$
is zero.
Now we consider a nontrivial extension $F$ in $\Ext^1(k(P),F_P)$
$$(\Scal_P) \qquad 0 \to F_P \to F \rarpa{\pi_P} k(P) \to 0 \,.$$
The scheme $\pdop(F)$ parameterizes surjections $\pi:F \to k(Q)$ from $F$ to torsion
sheaves of length one.
The subset of $\pdop(F)$ where $H^*(X, \ker(\pi) \otimes \Ecal_{c_0})=0$ is open
and not empty because
it contains $\pi_P$.
Thus, for a general point $Q$ of $X$, there exists a short exact sequence
$$(\Scal_Q) \qquad 0 \to F_Q \to F \to k(Q) \to 0$$
with $H^*(X,F_Q \otimes \Ecal_{c_0})=0$.
Applying the functor $R^* q_* (\Ecal_{r,i} \otimes  p^*(-)) $ to the short exact sequences
$(\Scal_P)$ and $(\Scal_Q)$, we obtain that $q_*(p^*F \otimes \Ecal)$ is isomorphic to
$\Ecal_{r,i}|_{C_{r,i} \times \{ P \}}$ and to
$\Ecal_{r,i}|_{C_{r,i} \times \{ Q \}}$ as well.

\tip
{\em Step 6: The filtration $F^*(\Ecal_{r,i})$ on the vector bundle $\Ecal_{r,i}$.}\\
If we consider the Harder--Narasimhan filtration on $G:=q_*(p^*F \otimes \Ecal)$,
then the graded summands need not be to be simple bundles. We consider a slight generalization
by taking $G_1$ to be a subsheaf of $G$ which is stable of maximal possible slope.
Defining $F^1(\Ecal_{r,i}) := G_1 \boxtimes q_*\shom(p^*G_1,\Ecal_{r,i})$, and 
$F^l(\Ecal_{r,i}) := \pi_1^{-1}( F^{l-1}(\Ecal_{r,i}/ F^1(\Ecal_{r,i})))$,
where $\pi_1$ is the surjection from $\Ecal_{r,i} \to \Ecal_{r,i}/ F^1(\Ecal_{r,i})$, we obtain
a filtration
$0 = F^0(\Ecal_{r,i}) \subset F^1(\Ecal_{r,i}) \subset \ldots \subset  F^k(\Ecal_{r,i}) = \Ecal_{r,i}$ 
on $\Ecal_{r,i}$ with the property that the $j$th graded object
$\gr^j(\Ecal_{r,i}) := F^j(\Ecal_{r,i}) / F^{j-1}(\Ecal_{r,i})$ is of the form $G_j \boxtimes F_j$.
By definition the slopes $\mu_j:=\mu(G_j)=\frac{\deg( G_j)}{\rk(G_j)}$ form a decreasing
sequence $\mu_1 \geq \mu_2 \geq \ldots \geq \mu_k$.

Restricting the filtration $F^*(\Ecal_{r,i})$ to a fiber $p^{-1}(x)$ of $p$, we obtain
a filtration of $G$ which does not depend on the choice of $x \in X$.
The restricted vector bundle $F^j(\Ecal_{r,i})|_{p^{-1}(x)}$ appears in the Harder--Narasimhan
filtration of $G$, if and only if $\mu_j > \mu_{j+1}$.
Therefore, we use $\HNF^*(\Ecal_{r,i})$ to name the subfiltration of the filtration
$F^*(\Ecal_{r,i})$ consisting of those $F^j(\Ecal_{r,i})$ with $\mu_j > \mu_{j+1}$.

\tip
{\em Step 7: Numerical invariants of the filtration  $F^*(\Ecal_{r,i})$.}\\ 
In the Grothendieck group $K(C_{r,i} \times X)$, we can identify $\Ecal_{r,i}$ with
the direct sum of the graded objects $\gr^j(\Ecal_{r,i})$:
$$[\Ecal_{r,i}] = \sum_{j=1}^k [\gr^j(\Ecal_{r,i})] = \sum_{j=1}^k [ G_j \boxtimes F_j ]\,.$$
Since the Chern character of the product $G_j \boxtimes F_j$ is given by
$$\begin{array}{rcl}
\ch(G_j \boxtimes F_j)&=& q^*\ch(G_j) . p^* \ch(F_j) \\
&=& \rk(G_j)\cdot \rk(F_j) + \left[ \rk(G_j) p^* c_1(F_j) + \rk(F_j) q^* c_1(G_j) \right]
+ p^* c_1(F_j). q^*c_1(G_j) \,,
\end{array}$$
we deduce the equality
$\int_{C_{r,i} \times X} \ch(\Ecal_{r,i})
= \sum_{j=1}^k \deg(G_j) \cdot \deg(F_j)$.
In step 3, we identified the left hand side with $\frac{-1}{2r} \deg(L_1)$.
Thus, we have
\begin{equation}
\sum_{j=1}^k \deg(G_j) \cdot \deg(F_j)  =  0
\end{equation}
The degree $\deg_X(F^j(\Ecal_{r,i})|_{q^{-1}(c)})$
of $F^j(\Ecal_{r,i})$, restricted to a fiber of $q$, is given by
$\deg_X(F^j(\Ecal_{r,i})|_{q^{-1}(c)})=\sum_{m=1}^j \rk(G_m) \cdot \deg(F_m)$.
Since the restriction of $\Ecal_{r,i}$ to a fiber of $q$ is semistable
of degree zero, we deduce that
\begin{equation}
A_j:=\sum_{m=1}^j \rk(G_m) \cdot \deg(F_m) \leq 0\,,
\end{equation}
and $A_k=0$.
Having in mind that $\mu_j:=\frac{\deg(G_j)}{\rk(G_j)}$, we rewrite equation (1)
$$0 = \sum_{j=1}^k \mu_j \cdot \rk(G_j) \cdot \deg(F_j)
= \sum_{j=1}^k A_j \cdot (\mu_j-\mu_{j+1}) \,,$$
where we set $\mu_{k+1}=0$.
The inequalities (2) and $\mu_j \geq \mu_{j+1}$ imply, therefore,
that $A_j=0$ whenever $\mu_j > \mu_{j+1}$.

\tip
{\em Step 8: $S$--equivalence of all bundles parametrized by $C_{r,i}$.}\\
The conclusion of the preceding steps is, that each quotient\\
$\gr_\HNF^j:= \HNF^j(\Ecal_{r,i})/\HNF^{j-1}(\Ecal_{r,i})$ is
on the one hand semistable of degree zero
when restricted to the fibers of $q$.
On the other hand, the semistable vector bundle on
$C_{r,i}$ which we obtain by
restricting $\gr_\HNF^j$ to a fiber $p^{-1}(x)$ of $p$ does not depend
on the chosen $x \in X$. Now the situation is symmetric in the sense
that $\gr_\HNF^j$ is semistable on all fibers of $p$ and $q$.
Thus, analogously to step 5, we deduce that the restrictions
$\gr_\HNF^j|_{\{ c_1 \} \times X }$ and $\gr_\HNF^j|_{\{ c_2 \} \times X }$
are isomorphic for all geometric points $c_1,c_2 \in X$.
In other terms, the direct sum $\oplus_{j=1}^l \gr_\HNF^j$
of the graded objects gives the same direct
sum of semistable vector bundles of degree zero on each fiber of $q$.
In short: All vector bundles parametrized by $C_{r,i}$ are $S$--equivalent.
\qed

\neu{Prop3}{\bf Theorem.} {\em
Let $\iota: Y \to X$ be a morphism of a smooth curve $Y$ to $X$.
We obtain a morphism $\iota_{M_{r,i}} : \tilde M_{r,i}(X) \to M_r(Y)$
by the pull back of vector bundles.
The following two conditions are equivalent:

\begin{tabular}{rl}
\hspace{2em} (1) & The degree of $\iota^*\Lcal_{r,i}$ is zero;\\
(2) & The morphism $\iota_{M_{r,i}}$ maps $\tilde M_{r,i}(X)$ to a point.
\end{tabular}}

\proof
We consider the morphisms
$\xymatrix{C_{r,i} \ar@{^(->}[r]^-{\alpha} & \tilde M_{r,i}(X)
\ar[r]^{\iota_{M_{r,i}}} & M_r(Y)}$. Theorem \ref{PROP}.\ref{Prop2}
implies that the degree of $\iota^*\Lcal_{r,i}$ is zero, if and only if
$\iota_{M_{r,i}}(C_{r,i})$ is a point. This implies the theorem
because $\tilde M_{r,i}(X)$ is irreducible, and $C_{r,i}$ is
the intersection of ample divisors.
\qed

\neu{Prop4}{\bf Corollary.} {\em
The equivalence classes of the nef line bundles $\Lcal_{r,i}$ and $\Lcal_r$
with respect to $\sim$ (see \ref{PROP}.\ref{Prop1})
neither depend on the choice of $C_{r,i} \subset \tilde M_{r,i}$
nor on the choice of the vector bundle $\Ecal_{r,i}$ on $C_{r,i} \times X$.
Furthermore, these equivalence classes
are independent of the chosen polarization $H$ on $X$.
}
\qed

\neu{Prop5}{\bf Corollary.} {\em
The line bundle $\Lcal_r$ on $X$ is nef.
For a morphism $\iota: Y \to X$ of a smooth curve $Y$ to $X$, we have
$\deg(\iota^* \Lcal_r)=0$, if and only if the morphism
$M_r(X) \to M_r(Y)$ is locally constant.}
\qed

\neu{Prop6}{\bf Proposition.} {\em
The line bundles $\Lcal_r$ satisfy the inequality
$\Lcal_r \preceq \Lcal_{r_2}$, for $r_1 \leq r_2$.\\
There exists a number $R \in \ndop$ such that $\Lcal_r \preceq \Lcal_R$ for all $r$.
}

\tip
\proof If $r_1 < r_2$, then we have an embedding of moduli spaces
$M_{r_1} \to M_{r_2}$ given by $[E] \mapsto [E \oplus \Ocal_X^{\oplus (r_2-r_1)}]$.
Thus, we deduce from theorem \ref{PROP}.\ref{Prop3}
that the inequality $\Lcal_{r_1} \preceq \Lcal_{r_2}$ holds.
We have seen in \ref{PROP}.\ref{Prop1} that in the chain
$\Lcal_1 \preceq \Lcal_2 \preceq \ldots \preceq \Lcal_r \preceq \ldots$ there are
at most $\rho(X)$ strict inclusions.
This proves the second assertion of the proposition.
\qed

\neu{Prop7}{\bf The line bundle $\Lcal_\infty$.}
We use the name $\Lcal_\infty$ for the line bundle $\Lcal_R$
of the above proposition. When referring to this line
bundle, we should be aware that $\Lcal_\infty$ is only
a class in \{nef line bundles\}$/ \sim$.
Considering these equivalence classes (with the discrete topology),
we have $\lim_{r \to \infty} \Lcal_r \sim \Lcal_\infty$.

In the following theorem, we have summarized the results of this section.

\neu{Prop8}{\bf Theorem (Properties of the line bundles $\Lcal_r$)}
{\em 
Let $X$ be a projective variety. We have an infinite sequence of nef line bundles
$\Lcal_1 \preceq \Lcal_2 \preceq \ldots \preceq \Lcal_r \preceq \ldots$
and a nef limit line bundle $\Lcal_\infty$ with
$\lim_{r \to \infty} \Lcal_r \sim \Lcal_\infty$
on $X$ such that for any morphism $\iota:Y \to X$ of a smooth curve $Y$ to $X$
the following conditions are equivalent:

\tip
\begin{tabular}{rp{13cm}}
\hspace{1em} (1) & $\deg(\iota^* \Lcal_r)=0$; \\
(2) & The restriction morphism $M_r(X) \to M_r(Y)$ is locally constant;\\
(3) & For any connected scheme $Z$ and every vector bundle $\Ecal$ on $Z \times X$
parameterizing semistable rank $r$ vector bundles on $X$ with trivial Chern classes,
the pull back $(\id_Z \times \iota)^*\Ecal$
parameterizes only one $S$--equivalence class on $Y$;\\
(4) & Modulo conjugation only finitely many $\U(r)$ representations of $\pi_1(Y)$ are
induced by those of $\pi_1(X)$.\\
\end{tabular}

}

\section{The generalized Albanese morphisms}\label{GA}
\neu{Ga1}
{\bf The construction of the generalized Albanese morphism.}
If $\Lcal_r$ or some power of it were base point free,
then it would define a morphism $\psi: X \to \pdop^m$.
Let $\varphi: X \to \Alb_r$ be the Stein factorization of
$\psi$, i.e.,
$\varphi$ is surjective with connected fibers.
Two geometric points $x$ and $x'$ of $X$ have the same image under $\varphi$,
if and only if $x \sim_{\Lcal_r} x'$.
By Theorem \ref{PROP}.\ref{Prop8}, the map $\varphi$ would
meet the requirements of a generalized Albanese variety.
Indeed, a curve $\xymatrix{Y \ar[r]^\iota &X}$ would be contracted by $\varphi$,
if $\deg_Y(\iota^* \Lcal_r)=0$. This means (by theorem \ref{PROP}.\ref{Prop8})
that all families of semistable rank $r$ vector bundles on $X$ with trivial
Chern classes become constant 
when restricted to $Y$, or only finitely many representation classes modulo
conjugation of $\pi_1(Y)$ in
$\U(r)$ are induced by representations of $\pi_1(X)$.

If no line bundle $L_r$ with $L_r \sim \Lcal_r$ is base point free
(Note, that $L^{\otimes k} \sim L$, for all $k >0$.),
then Tsuji's nef reduction theorem provides us
with a rational version of the generalized Albanese variety up
to birational equivalence.
In this case we obtain only a birational model of the Albanese
morphism and variety.

\neu{Ga2}
{\bf Theorem.} (see Theorem 2.1 in \cite{the8},
see also \cite{Tsuji})
{\em
There exists a dominant rational map
$\xymatrix{X \ar@{-->}[r]^-{\alb_r} & \Alb_r(X)}$ with connected fibers such that:

\begin{tabular}{rl}
\hspace{2em} (1) & The line bundle $\Lcal_r$ is numerically trivial on all
compact fibers $F$ of $\alb_r$\\
& of dimension $\dim(X)-\dim(\Alb_r(X))$;\\
(2) & For every general point $x \in X$ and every
irreducible curve $C$ passing\\
&through $x$ with $\dim(\alb_r(C))>0$, we have $C.\Lcal_r >0$;\\
(3) & There exist compact fibers of $\alb_r$.\\
\end{tabular}

Furthermore, the pair $(\alb_r,\Alb_r(X))$ is uniquely
determined up to birational equivalence. 
}

\neu{Ga3}{\bf The chain of generalized Albanese morphisms.}
Even though we end up with an infinite sequence of rational
morphisms $\xymatrix{X \ar@{-->}[r]^-{\alb_r} & \Alb_r(X)}$, for each $r \in \ndop$,
there are at most $\rho(X)$ different generalized Albanese morphisms,
since for almost all $r \in \ndop$, we have $\Lcal_r \sim \Lcal_{r+1}$.
Since $\Lcal_r \preceq \Lcal_{r+1}$, we get a rational morphism
$\xymatrix{\Alb_{r+1}(X) \ar@{-->}[r] & \Alb_r(X)}$.
So, we end up with the following commutative diagram:
$$\xymatrix{& \Alb_\infty(X) \ar@{-->}[d] \ar@{-->}[dr] &
X \ar@{-->}[l]_-{\alb_\infty} \ar@{-->}[d]^-{\alb_r} \ar@{-->}[dl]
\ar[drr]^-{\alb_1}\\
\stackrel{\mbox{ }}{\ldots} \quad \ar@{-->}[r] &
\Alb_{r+1}(X) \ar@{-->}[r] & \Alb_r(X) \ar@{-->}[r]
& \quad \stackrel{\mbox{ }}{\ldots} \quad \ar@{-->}[r] & \Alb_1(X)}$$

\neu{Ga4}{\bf Proposition. (Functoriality)}
{\em 
If $\psi: X \to X'$ is a morphism of projective varieties,
then we have $\psi^* \Lcal'_r \preceq \Lcal_r$
for all $r \in \ndop \cup \{ \infty \}$. Therefore,
we have commutative diagrams
$$\xymatrix{X \ar[r]^\psi \ar@{-->}[d]^-{\alb_r} & Y \ar@{-->}[d]^-{\alb_r} \\
\Alb_r(X) \ar@{-->}[r] & \Alb_r(Y).} $$
}

\proof
Suppose that $\iota:Y \to X$ is a morphism from a smooth curve to $X$
with $\deg(\iota^* \psi^* \Lcal'_r) >0$. This implies by theorem
\ref{PROP}.\ref{Prop3} that there exists a family $\Ecal'$ of rank $r$ vector
bundles with trivial Chern classes on $X'$
parametrized by a connected scheme $Z$ such that the pull back
$(\id_Z \times (\psi \circ \iota))^*\Ecal'$ of $\Ecal'$ to $Z \times Y$
parameterizes different $S$--equivalence classes on $Y$.  
However, then the family $\Ecal=(\id_Z \times \psi)^*$ also parameterizes
rank $r$ vector bundles on $X$ which pull back to non $S$--equivalent classes.
Consequently, again by theorem \ref{PROP}.\ref{Prop3}, $\deg(\iota^* \Lcal_r) >0$.
\qed

\neu{Ga5}{\bf Proposition.}
{\em
If $\psi:X \to X'$ is an \'etale morphism of projective varieties,
then we have $\Lcal_r \preceq \psi^* \Lcal'_{r \cdot \deg(\psi)}$
for all $r \in \ndop$.
Furthermore, $\Lcal_\infty \sim \psi^* \Lcal'_\infty$.
}

\tip
\proof
The proof is analogous to the preceding one. We simply have to consider
the push forward of a family of rank $r$ vector bundles on $X$
to $X'$. The statement about $\Lcal_\infty$ is easily obtained from
$\Lcal_r \preceq \psi^* \Lcal'_{r \cdot \deg(\psi)}
\preceq \Lcal_{r \cdot \deg(\psi)}$ and the fact that $\Lcal_r \sim \Lcal_\infty$
for $r \gg 0$.

\neu{Ga6}{\bf Remark.}
Since the first Albanese variety and morphism exist,
it would be enough to have the base point freeness of
some power of $\Lcal_r$ on every fiber of the Albanese
variety to obtain a regular morphism
$$\xymatrix{ X \ar[rr]^-{\alb_r} \ar[drr]_-{\alb} && \Alb_r(X)\ar[d]\\&&\Alb(X)\,.}$$
This is used to study the
generalized Albanese morphisms in the case of algebraic surfaces
in the next section.

\section{The case of algebraic surfaces}\label{SURFACE}

\neu{sf0}
We consider here the case of a polarized projective
algebraic surface $(X,H)$.
In this case, we can make a much stronger statement
than theorem \ref{GA}.\ref{Ga2}.
The generalized Albanese morphism is well understood 
if $X$ is not of general type
(see \ref{SURFACE}.\ref{sf6}).
The rest of this section
(\ref{SURFACE}.\ref{sf1} -- \ref{SURFACE}.\ref{sf5}) is devoted to
the proof of the

\neu{sfmain}{\bf Theorem. \/}{\em
Let $(X,H)$ be a polarized projective surface.
Then there exists a surjective morphism
$\alb_r:X \to \Alb_r(X)$ with connected fibers,
and an effective divisor $D$ on $X$,
such that for all morphisms $\iota:C \to X$
of irreducible curves
with $\iota(C) \not \subset D$ the following conditions
are equivalent.\\
\begin{tabular}{rl}
\hspace{2em} (1) & $\alb_r(\iota(C))$ is a point;\\
(2) & The associated morphism
$\xymatrix{\Hom(\pi_1(X),\U(r)) \ar[r]^-{\iota^*} &
\Hom(\pi_1(C),\U(r))}$\\
&modulo conjugation has a finite image;\\
(3) & For any base scheme $S$ and any rank $r$
vector bundle $E$ on $X \times S$\\
&such that for each $s$ $E_s$ is semistable with
numerically trivial Chern\\
& classes, the pull back of $E$ to $C \times S$ is
a family of S--equivalent vector\\
&bundles.
\end{tabular}

The divisor $D$ of exceptions can be written in the form
$D = C_1 + C_2 +\ldots + C_l$,
where the $C_i$ are irreducible and form a basis
of a proper subspace of the rational N\'eron--Severi
vector space $NS(X) \otimes \qdop$.
In particular, we have
$l< \dim_\qdop(NS(X) \otimes \qdop)$.}

\neu{sf1}
{\bf Preparations for the proof.}
We consider the nef line bundle $\Lcal_r$ on $X$ satisfying the equivalence of
theorem \ref{PROP}.\ref{Prop3}. 
There are two extreme cases where the proof is a simple remark:
When $C.\Lcal_r > 0$, for all curves $C$, then we set $\alb_r$ to be the identity
morphism of $X$.
If $C.\Lcal_r=0$ for all curves, then we set $\Alb_r(X)=\Spec(\cdop)$, and we are
finished.

Thus, we assume from now on that $\Lcal_r $ is a numerically nontrivial nef line bundle
vanishing on a nonempty set $\{ C_i \}_{i \in I}$ of irreducible curves.
It follows from the construction that all the curves $C_i$ in this family
are contracted to points by the classical Albanese morphism.

$\Lcal_r$ being a nef divisor is the limit of ample
divisors (with rational coefficients) which yields $\Lcal_r^2 \geq 0$.
This implies the following lemma.

\neu{sf2}
{\bf Lemma.}
{\em If $C \subset X$ is an effective divisor with $C^2>0$, then $C.\Lcal_r >0$.}

\tip
\proof We assume the contrary.
Since $\Lcal_r$ is nef this means that $C.\Lcal_r=0$.
The Hodge index theorem (see IV.2.15 in \cite{BPV}) implies that $\Lcal_r^2 \leq 0$
with equality only when $\Lcal_r$ is torsion. From $\Lcal_r^2 \geq 0$ and the assumption that
$\Lcal_r$ is numerically nontrivial we derive a contradiction.
\qed

\tip
The proof of theorem \ref{SURFACE}.\ref{sfmain}
is subdivided into three cases depending on
the dimension of $X$ in the Albanese variety.
This is by definition the dimension of $\Alb_1(X)$.

\neu{sf3}
{\bf The nef reduction for $\Alb_1(X) = \Spec(\cdop)$.}
We consider the curves $\{ C_i \}_{i \in I}$ as vectors in the rational
Ner\'on--Severi space $\NS_\qdop(X)$.
The dimension of this vector space is the Picard number $\rho(X)$ of $X$.
The points $\{ C_i \}_{i \in I}$ lie on the hyperplane 
$\{ C \in \NS_\qdop(X) \, | \, C. \Lcal_r=0\}$.
Suppose there would be a nontrivial linear relation among the $C_i$.
If the number of these curves exceeds $\rho(X)-1$,
then we have at least one such relation.
We write the linear relation
$a_1C_1+\ldots a_mC_m = a_{m+1}C_{m+1}+ \ldots a_MC_M$ with positive rational $a_i$,
and $C_i$ different from $C_j$ whenever $i \not=j$.

After multiplication with a positive integer,
we may assume the $a_i$ to be integers.
We set $D_1:=a_1C_1+\ldots a_mC_m$ and $D_2:=a_{m+1}C_{m+1}+ \ldots a_MC_M$.
Since $D_1$ and $D_2$ coincide in $\NS_\qdop(X)$,
their difference is torsion in
the Ner\'on--Severi group.
So again, after multiplying with an integer,
we may assume that the effective Cartier divisor classes $D_1$ and $D_2$ coincide.
(Here we use the fact that the
Picard torus, the dual of the Albanese torus, is trivial.)

$D_1^2 =D_1.D_2 \geq 0$, because $D_1$ and $D_2$ have no common components.
In view of lemma \ref{SURFACE}.\ref{sf2}, and $D_1.\Lcal_r=0$,
we conclude that $D_1^2=0$.
This implies that $D_1$ and $D_2$ are disjoint.

Consequently, the line bundle $L:=\Ocal_X(D_1) \cong \Ocal_X(D_2)$ has two linearly
independent sections which do not intersect.
Thus, $L$ is base point free
and defines a morphism whose Stein factorization we denote by $\alb_r: X \to \Alb_r(X)$. 

On the one hand,
we have that $F.\Lcal_r=0$ for all fibers of $\alb_r$.
On the other hand, suppose $C.\Lcal_r=0$ for a curve $C \subset X$.
Let $F$ be an irreducible fiber of $\alb_r(X)$.
If $C$ were not contained in a fiber, then we would have $(C+mF)^2>0$ for $m \gg0$.
However, we have $(C+mF).\Lcal_r=0$ which contradicts lemma  \ref{SURFACE}.\ref{sf2}.
Thus, each curve $C$ with $C.\Lcal_r=0$ is contained in a fiber.

This means that the effective divisor $D$ of \ref{SURFACE}.\ref{sfmain}
can be taken to be the empty set once we have a linear relation between the
$\{C_i \}_{i \in I}$ in $\NS_\qdop(X)$. Since the resulting morphism $\alb_r$
contracts all these curves, we conclude that $\alb_r$ does not depend on 
the chosen linear relation.

\neu{sf4}
{\bf The nef reduction when $\Alb_1(X)$ is a curve.}
We consider the morphism $\alb_1:X \to \Alb_1(X)$.
This morphism is the Stein factorization of the classical Albanese morphism.
It follows from $\Lcal_1 \preceq \Lcal_r$,
that each curve $C$ with $C.\Lcal_r = 0$ is contained in a fiber
of this morphism.
Let $F$ be the generic fiber of $\alb_1$.

If $F.\Lcal_r=0$, then all curves $C$ with $C.\Lcal_r=0$ are contracted by $\alb_1$.
Consequently we set $\alb_r=\alb_1$ and theorem \ref{SURFACE}.\ref{sfmain} is proven.

We suppose now that $F.\Lcal_r$ is positive.
The set of curves $\{C_i\}_{i \in I}$ consists of components of reducible fibers
of the morphism $\alb_1$.
We will show that this set is not only finite but linearly independent in $\NS_\qdop(X)$.
Indeed, if there were a linear relation, then we would obtain 
(see \ref{SURFACE}.\ref{sf3}) an effective divisor $D_1=a_1C_1 + \ldots a_mC_m$.
This divisor satisfies $D_1^2=0$, $D_1.\Lcal_r=0$,  and consists of fiber components.
This contradicts Zariski's lemma (see Lemma III.8.2 in \cite{BPV}),
because $F.\Lcal_r>0$ for all fibers $F$.
This shows that theorem \ref{SURFACE}.\ref{sfmain} holds when setting $\alb_r = \id_X$.
 
\neu{sf4a}
{\bf Remark.}
The finite collection of curves $\{C_i\}_{i \in I}$ must have a negative definite
intersection matrix, because of Zariski's lemma.
Thus, by Grauert's criterion (Theorem III.2.1 in \cite{BPV}),
there exists a contraction of these curves.
However, this contraction is not necessarily a projective morphism.
If it were, we could take this contraction to be our generalized
Albanese morphism $\alb_r$.

\neu{sf5}
{\bf The nef reduction for $\dim(\Alb_1(X))=2$.}
In this case, the generalized Albanese morphism $\alb_1:X \to \Alb_1(X)$
contracts finitely many curves. Among those curves are the $\{C_i\}_{i \in I}$
which are numerically trivial with respect to $\Lcal_r$.
The intersection matrix of these  $\{C_i\}_{i \in I}$ is negative definite.
This yields that these curves form a basis of a proper subspace of $\NS_\qdop(X)$.
Again, just setting $\alb_r=\id_X$, the assertions of \ref{SURFACE}.\ref{sfmain}
are fulfilled.
As before remark \ref{SURFACE}.\ref{sf4a} applies.

\neu{sf6}{\bf Surfaces of Kodaira dimension less than 2.}
Let $X$ be a projective algebraic surface of Kodaira
dimension $\kappa(X) \leq 1$.
We assume that $X$ is minimal.
This is not a restriction, because the fundamental group of rational curves is zero.
The next table gives the generalized Albanese morphism
for these surfaces following the Enriques-Kodaira classification (see VI. in \cite{BPV}).
The last two rows are perhaps the most interesting ones.
They show that the generalized Albanese morphisms may reveal
more of the surface than the classical one.
We assume in these two rows that the surface $X$ is not a
product of two curves.

\begin{center}
\begin{tabular}{l|l|l}
$\kappa(X)$ & class of $X$ & the generalized Albanese morphism\\
\hline
$-\infty$ & rational surfaces & $X \to \Spec(\cdop)$ \\
\hline
$-\infty$ & ruled surfaces & $X \to B$\\
& $X \to B$ with $g(B) \geq 1$\\
\hline
0& Enriques surfaces & $X \to \Spec(\cdop)$\\
\hline
0& K3 surfaces & $X \to \Spec(\cdop)$\\
\hline
0 & tori & $X \to X$\\
\hline
0 & hyperelliptic surfaces &
$\Alb_1(X)$ is an elliptic curve, whereas \\
&&$\Alb_r(X) \cong X$, for $r>1$, see also \S \ref{EXA}.\\
\hline
1 & properly elliptic surfaces & $\Alb_1(X)$ is an algebraic curve, and \\
&&$\Alb_r(X) \cong X$, for $r>1$.\\
\end{tabular}
\end{center}

\section{An example}\label{EXA}
\neu{exa1}{\bf The group $G$.}
We consider the group $G$ which as a set is $\zdop^4$ with group structure
given by
$(a,b,c,d) (a',b',c',d')=(a+a',b+b',c+(-1)^ac',d+(-1)^ad')$.
The group $G$ has the four generators $g_0=(1,0,0,0)$,
$g_1=(0,1,0,0)$,
$g_2=(0,0,1,0)$,
and $g_3=(0,0,0,1)$ satisfying the six relations:
$$g_0g_1=g_1g_0 \quad
g_0g_2=g_2^{-1}g_0 \quad g_0g_3=g_3^{-1}g_0 \quad g_1g_2=g_2g_1
\quad g_1g_3=g_3g_1 \quad g_2g_3=g_3g_2\,.$$
$G$ may be described as the semidirect product
$(\zdop^2) \ltimes (\zdop^2)$ where the homomorphism
$(\zdop^2) \to \Aut(\zdop^2)$ is given by
$(a,b) \mapsto (-1)^a \id_{\zdop^2}$.

\neu{exa2}{\bf Characters of $G$.}
It is easy to check that the commutator subgroup $[G,G]$ is generated by
$g_2^2$ and $g_3^2$.
Thus,
we obtain $G/[G,G] \cong \zdop^2 \oplus (\zdop/2\zdop)^2$.
Consequently,
each group homomorphism $\tau:G \to \U(1)$ must send $g_2$
and $g_3$ to $\pm 1$.
Thus,
in a continuous family $\tau_t$ of such morphisms the images of
$g_2$ and $g_3$ are locally constant.
We conclude that there are four  families of $\U(1)$-representations of
$G$.
One is given by
$$\tau_{t_1,t_2}(g_0)= \exp(t_1 \cdot i) \quad
\tau_{t_1,t_2}(g_1)= \exp(t_2 \cdot i) \quad
\tau_{t_1,t_2}(g_2) \equiv 1 \equiv \tau_{t_1,t_2}(g_3)$$
with $(t_1,t_2) \in (\rdop/2 \pi \zdop)^2$.
The other three representations of $G$ in $\U(1)$ are obtained by changing
the image of $g_2$ or $g_3$ to $-1$.

\neu{exa4}{\bf $G$ is the fundamental group of a hyperelliptic surface.}
We consider two elliptic curves $(E_1,0)$ and $(E_2,0)$ with a nonzero
two-torsion point $e_1 \in E_1$.
The fix point free action of $\zdop/2\zdop$ on $E_1 \times E_2$ generated
by
$(z_1,z_2) \mapsto (z_1+e_1,-z_2)$ induces a smooth projective quotient
$X$,
and a $2:1$ covering $p:E_1 \times E_2 \to X$.
This $X$ is a hyperelliptic surface (see \cite{BPV} V.5).
The fundamental group $\pi_1(X)$ is isomorphic to $G$ and the commutative
subgroup $2\zdop \times \zdop^3$ of index two corresponds to the cover $p$.
Denoting the $\zdop/2\zdop$-quotient of $E_1$ by $e_1$ by $\bar E_1$, we
obtain a
morphism $\alpha:X \to \bar E_1$.
The induced map on the first homology groups
$$\alpha_*:\,\left[H_1(X,\zdop) \cong \zdop^2 \oplus (\zdop/2\zdop)^2 
\right]
\longrightarrow
\left[H_1(\bar E_1,\zdop) \cong \zdop^2 \right]$$
is an isomorphism between the free parts whereas the torsion part of
$H_1(X,\zdop)$ is mapped to zero.
Thus,
the morphism $\alpha$ is the Albanese morphism for the surface $X$.

\neu{exa6}{\bf The generalized Albanese morphism $\alb_2(X)$ differs from
the ordinary Albanese morphism.}
Let $\Lcal_2$ be the line bundle defining the second generalized Albanese morphism
$\alb_2(X)$ on $X$. And $\tilde \Lcal_1$ be the corresponding line bundle on $E_1 \times E_2$.
Since $\Alb(E_1 \times E_2)= E_1 \times E_2$, the line bundle $\tilde \Lcal_1$ is ample.
From proposition \ref{GA}.\ref{Ga5} we deduce that $\tilde \Lcal_1 \preceq p^* \Lcal_2$.
Thus $\Lcal_2.C >0$ for all curves $C$ in $X$.

Since $\Lcal_2 \preceq \Lcal_r$ for all $r \geq 2$,
the equivalence relation $\sim_{\Lcal_r}$ defined by the line bundles $\Lcal_r$
is the identity relation.
Thus, $\alb_r = \id_X: X \to X $ for all $r \geq 2$.
On the other hand, the classical Albanese variety has only dimension one.

\neu{exa7}{\bf A class of examples.}
The above example is a special case of the following class of examples.
Let $G$ be a finite group with $|G|$ elements and $C_1$, $C_2$ be two
smooth projective curves with a $G$ action such that 

\tip
\begin{tabular}{rp{12cm}}
\hspace{2em} (1) & The genera $g_{C_1}$ and $g_{C_2}$ are positive;\\
(2) & $G$ acts free on $C_1$, i.e the quotient map $C_1 \to C_1/G$ is \'etale; \\
(3) & There are no $G$--invariant global sections in $H^0(C_2,\omega_{C_2})$.
This is equivalent to $C_2/G \cong \pdop^1$.
\end{tabular}

\tip
We obtain a free $G$--action on $C_1 \times C_2$.
Let $X:=(C_1 \times C_2)/G$ be the quotient of this action
and $p:C_1 \times C_2 \to X$ the projection.
Since $C_1 \times C_2$ is embedded into its Albanese variety,
we deduce from proposition \ref{GA}.\ref{Ga5}, that $\Lcal_{|G|}.C >0$
for all curves $C \subset X$. Thus, $\alb_r = \id_X: X \to X$ for all
$r \geq |G|$, whereas the classical Albanese morphism is just the map
to the curve $C_1/G$.

{\small
Georg Hein\\
Freie Universit\"at Berlin,
Institut f\"ur Mathematik II, Arnimallee 3,
D-14195 Berlin, Germany\\
{\tt ghein@math.fu-berlin.de}}
\end{document}